\documentclass{article}
\usepackage{latexsym}
\usepackage{amsthm}  
\usepackage[all]{xy}
\usepackage{stmaryrd}
\usepackage{amssymb}
\usepackage{mathptmx}

\newcommand{\mbb}[1]{\mathbf{#1}}

\newcommand{\msc}[1]{\mathcal{#1}}
\newcommand{\mrm}[1]{\mathrm{#1}}

\newcommand{\Z}{\mbb{Z}}
\newcommand{\Q}{\mbb{Q}}
\newcommand{\C}{\mbb{C}}
\newcommand{\Pro}{\mbb{P}}

\newcommand{\sh}[1]{{\msc{#1}}}
\newcommand{\fn}[3]{{#1:#2\rightarrow #3}}

\newcommand{\Proj}{\mrm{Proj~}}

\theoremstyle{definition}
\newtheorem{definition}{Definition}[section]
\newtheorem{remark}[definition]{Remark}

\newtheorem{example}[definition]{Example}
\theoremstyle{plain}
\newtheorem{theorem}[definition]{Theorem}

\newtheorem{lemma}[definition]{Lemma}
\newtheorem{proposition}[definition]{Proposition}

\date{}		

\begin{document}

\title{Stable degenerations of symmetric squares of curves}
\author{Michael A. van Opstall}
\maketitle

\abstract{
The stable (in the sense of the relative minimal model program) degenerations
of symmetric squares of smooth curves of genus $g>2$ are computed. This 
information is
used to prove that the component of the moduli space of stable surfaces
parameterizing such surfaces is isomorphic to the moduli space of stable
curves of genus $g$.
}

\section{Introduction}

Suppose $C\rightarrow \Delta'$ is a family of smooth curves of genus greater
than two over a punctured
disk. Denote by $C^{(2)}_{\Delta'}$ the {\em fibered symmetric square} of 
this family, that is, the quotient of the fibered product $C\times_{\Delta'}
C$ by the $\Z_2$ action swapping the factors. The theory of moduli of stable
surfaces ensures that (after possibly taking a base change by a cover of
$\Delta$ totally ramified over 0) this family of surfaces may be completed
to a family of surfaces $S\rightarrow \Delta$ satisfying
\begin{enumerate}
\item $S_0$ has {\em semi-log canonical} (slc) singularities;
\item some reflexive power $\omega_{S/\Delta}^{[N]}$ of the relative
dualizing sheaf of the family is locally free;
\item $\omega_{S_0}^{[mN]}$ is an ample line bundle for some $m>0$.
\end{enumerate}
By analogy with the theory of moduli of stable curves, the process of finding
$S\rightarrow \Delta$ is called {\em stable reduction}. 

The main result of this article is a geometric description of these 
stable reductions. The first step is completing $C\rightarrow \Delta'$ to
a family of stable curves over $\Delta$ (base changes are ignored in this
non-technical discussion). The family $C_{\Delta}^{(2)}\rightarrow \Delta$ is 
not stable in general. A partial resolution of this family is given by the
Hilbert scheme $\mrm{Hilb}_2(C/\Delta)$ of length two subschemes in the
fibers. Geometrically, this replaces the points on $C_{\Delta}^{(2)}$ coming
from the quotient of a product of a node with itself by rational curves. 
The fibers of $\mrm{Hilb}_2(C/\Delta)$ have slc singularities, and the
relative dualizing sheaf is locally free, but the special fiber does not
yet in general have an ample dualizing sheaf. The last step of stabilization
is an exercise in geometry of products and symmetric squares of smooth 
curves.

This explicit description of the stable degenerations yields a description
of the irreducible component of the moduli space of stable surfaces which
is the closure of the moduli space of symmetric squares of curves of a 
given genus.

This article continues the program of studying moduli spaces of stable 
surfaces using compact moduli spaces of simpler objects to complete
families of surfaces. As in \cite{vo:mpc}, the moduli space of stable curves
is used here to ensure that a degeneration of smooth curves may be replaced
by a degeneration with special fiber at worst nodes. Moduli spaces of stable
curves with group action were used to study surfaces isogenous to a product
in \cite{vo:isog}. One could conceivably use the methods of Abramovich 
and Vistoli \cite{av:mf} to study stable degenerations of Kodaira
fibrations (see \cite{bpv:ccs2}, V.14). This method has been employed by
La Nave in the case of elliptic fibrations \cite{ln:mell}. Studying
symmetric squares of smooth curves by stable degeneration of curves has
yielded results on the cone of effective divisors of such surfaces; the
article \cite{ck:scgm} contains excellent pictures which may help the reader
in visualizing the constructions used here.

I thank Brendan Hassett for suggesting to me the interpretation of the
variety $\tilde{X}$ below as a Hilbert scheme.

\section{Stable surfaces and their moduli}

This section contains the relevant references and definitions from the
theory of moduli of stable surfaces which are necessary in what follows.
Since ampleness of the canonical class is more important for us than 
nonsingularity, we will consider families whose general member is a
{\em canonically polarized surface}, that is, a projective surface with
rational double points and ample dualizing sheaf. 

\begin{definition}
A {\em stable surface} is a connected, projective, reduced surface $S$ with 
semi-log canonical (slc) singularities such that $\omega_S^{[N]}$ (the
reflexive hull of the $N$-th power of the dualizing sheaf of $S$) is an
ample invertible sheaf for some $N$. 
\end{definition}

\begin{remark}
For the definition of slc singularities, see \cite{ksb:3f}, Definition 4.17,
or the equivalent, but differently formulated Definition 2.8 in \cite{a:mgnw}.
In the first article, all slc surface singularities are classified, and we will
use the classification more than the definition. For our purposes, the 
following information suffices:
\begin{enumerate}
\item Normal crossings singularities are slc.
\item Products of slc singularities are slc (in particular, the product of
nodal curves has slc singularities -- these are special cases of 
{\em degenerate cusps}).
\item Slc singularities are those which appear on the fibers of relative
canonical models of semistable reductions (see below).
\end{enumerate}
\end{remark}

It turns out that the naive definition of a family of stable surfaces --
a flat, proper morphism whose fibers are stable surfaces -- does not lead
to a separated moduli functor. One needs a further condition on the family.

\begin{definition}
A {\em family of stable surfaces} is a flat, proper morphism $\fn{f}{X}{B}$ 
whose fibers are stable surfaces and whose relative dualizing sheaf 
$\omega_{X/B}$ is $\Q$-Cartier, i.e., some reflexive power is a line bundle.
One also says that the morphism $f$ is $\Q$-Gorenstein.
\end{definition}

In \cite{k:pcm}, a stronger condition is required on a family of stable 
surfaces. For the purposes of this article, the weaker condition given here
suffices, since the weaker condition given here implies Koll\'ar's stronger
condition if the family has a smooth curve as base and canonically polarized
general fiber.

An essential fact for moduli theory of stable surfaces, which is proved using 
Mori theory is the following:

\begin{theorem}[Stable reduction]
Let $X\rightarrow B$ be a one-parameter family with a smooth base such that
$X_b$ is a canonically polarized surface for all $b\neq 0$ for some point 
$0\in B$. Then there exists
a finite map $B'\rightarrow B$ totally ramified over $0$ and a family of
stable surfaces $X'\rightarrow B'$ extending the pullback of the family
$X|_{B\backslash\{0\}}$ to $B'$. The special fiber of $X'$ is uniquely
determined by the original family.
\end{theorem}

\begin{proof}[Sketch of proof]
After a base change, the original family admits a semistable resolution,
so we may assume $X\rightarrow B$ is a family with smooth total space and
normal crossings divisors for fibers. In this case, a unique relative
canonical model for the morphism $X\rightarrow B$ exists (see Chapter 7 of
\cite{km:bgav} for a proof). This relative canonical model is a family of
stable surfaces. The uniqueness of the special fiber follows from the facts
that every pair of semistable resolutions is dominated by a third, and the
uniqueness of relative canonical models.
\end{proof}

Once one has established the existence of a coarse moduli space of finite
type (this depends on the boundedness theorem of Alexeev proved in
\cite{al:bk2}) for stable
surfaces (using the notion of family given above), this theorem implies that
the moduli space is proper (and in particular, separated). We will only
need the uniqueness part of the theorem, since our work will prove existence
for the special class of families which we consider.

A one-parameter family $X\rightarrow B$ with smooth base and canonically
polarized
general fiber as in the theorem will be called a {\em degeneration}. If
a degeneration $X\rightarrow B$ is a family of stable curves or stable 
surfaces, it will be called a {\em stable degeneration}.

\section{Punctual Hilbert schemes and Chow varieties of nodal curves}

Let $C\rightarrow B$ be a stable degeneration of curves. The
fibered symmetric square $C^{(2)}_B$ is isomorphic to the relative
Chow variety $\mrm{Chow}_{0,2}(C/B)$ of dimension zero, degree two cycles in 
the fibers of $C\rightarrow B$. Abbreviate $C^{(2)}_B$ by $X$.
The points of $X$ fall into
four classes:
\begin{enumerate}
\item Cycles consisting of two smooth points (possibly equal): such points are 
smooth on $X$ and on its fibers.
\item Cycles consisting of a smooth point and a node: such points are smooth
on $X$ and are normal crossings on its fibers.
\item Cycles consisting of two different nodes: such points are isolated
singular points on $X$ (analytically isomorphic to the vertex of a cone over 
a quadric surface) and degenerate cusps on its fibers.
\item Cycles consisting of the same node taken twice: such points are
isolated singular points on $X$ (analytically $\Z_2$ quotients of the cone 
over a 
quadric surface). Such singularities are analytically isomorphic to the
vertex of the cone over the cubic scroll. They are thus singular points
of the fibers as well.
\end{enumerate}

The cone over the cubic scroll is not $\Q$-Gorenstein, so the family
$X\rightarrow B$ cannot be $\Q$-Gorenstein if the original
family $C\rightarrow B$ contains nodal curves. The points of $X$ of
the fourth type considered above will be called the {\em bad points} of $X$. 

Let $\bar{X}\rightarrow X$ be the blowup of all bad points of $X$. It is
easy to check by computing in local coordinates that the fibers of
$f:\bar{X}\rightarrow B$ have only slc singularities (indeed, the only 
singular points in the fibers are the normal crossings and degenerate cusps 
described as points
of types 2 and 3 above). The stable reduction
of $X\rightarrow B$ will be the relative canonical model of this morphism,
that is,
\[
\Proj \bigoplus_{n=0}^\infty f_*\omega_{\bar{X}}^n.
\]

A first step towards a geometric description of this model is the relative
canonical model of $\bar{X}\rightarrow X$. Call this $\tilde{X}$. Example
2.7 of \cite{km:bgav} describes how to obtain $\tilde{X}$ directly from
$X$. This variety $\tilde{X}$ is a resolution of the singularities of $X$ 
with as small as possible exceptional locus, in this case, consisting of  
rational curves
rather than divisors. There is a preferred choice for such a resolution, which
is distinguished by having an ample relative canonical class. 

Let us recall the local description of $\tilde{X}$ from {\em loc. cit.}.
A bad point has a local analytic description as $xy-uv=0$ modulo the action
of $\Z_2$ taking  $x$ to $-x$ and $u$ to $-u$. In these coordinates, the 
parameter on the base is given by $xv-yu$, up to a constant factor. In these
local coordinates, we obtain $\tilde{X}$ by blowing up the ideal $(x,u)$ in
$\C[x,y,u,v]/(xy-uv)$, and then taking the quotient by the $\Z_2$-action
described above (which extends to the blowup). The resulting variety
has coordinate ring
\begin{eqnarray*}
\C[a_1,\ldots, a_7]/(a_5a_6-a_4a_7,~a_3a_6-a_1a_7,~a_2a_6-a_3a_7, \\
~a_3a_4-a_1a_5,~a_2a_4-a_3a_5,~a_1a_2-a_3^2),
\end{eqnarray*}
where the variables $a_6$ and $a_7$ are homogeneous. The central fiber of
this new family is obtained by setting $a_3=a_4a_5$.
One easily checks in coordinates that the
only singularities of the central fiber are normal crossings. So the
stable reduction improves the singularities of the central fiber. The total
space of the new family is smooth, and the canonical sheaf of
$\tilde{X}$ is ample relative to $\tilde{X}\rightarrow X$.

\subsection{Coordinate-free description of $\tilde{X}$}

There is another way to obtain this model, which is easier to describe.
However, for computations, the description given above it indispensible.
Let $C$ be a smooth curve. Then it is well-known that the Hilbert scheme
$\mrm{Hilb}_2(C)$ parameterizing length two subschemes of $C$ is isomorphic
to the Chow variety $\mrm{Chow}_{0,2}(C)$ of dimension zero, degree two
cycles. If $C$ has nodes, or more generally, plane curve singularities, 
each singularity supports a $\Pro^1$ (the projectivization of the tangent
space to this singularity) of length two structures, so the associated
Hilbert-to-Chow morphism can be described exactly as the morphism
$\tilde{X}\rightarrow X$ is described above. In fact, following our notation
above,

\begin{theorem}
\[
\tilde{X}\cong \mrm{Hilb}_2(C/B)
\]
\end{theorem}

\begin{proof}
By construction, the varieties $\tilde{X}$ and $\mrm{Hilb}_2(C/B)$ are
birational over $B$. Since $\tilde{X}$ is a relative canonical model,
there is actually a morphism $f:\mrm{Hilb}_2(C/B)\rightarrow \tilde{X}$,
which if not an isomorphism, must blow down the exceptional $\Pro^1$s in
the morphism $\mrm{Hilb}_2(C/B)\rightarrow X$.

However, the morphism $f$ is surjective, being projective. Consequently the
$\Pro^1$s occuring on $\mrm{Hilb}_2(C/S)$ are not blown down by this
morphism. 
\end{proof}

\subsection{A third description of $\tilde{X}$}

We will need a description of the special fiber $\tilde{X}_0$, or alternately,
$\mrm{Hilb}_2(C_0)$. In the fiber $X_0$, a sufficiently small analytic germ
at a bad point has three irreducible components (see Figure 1). One of these
components is distinguished in that it meets both of the others along curves;
the other pair of components meet only at the bad point. The local computation
given above shows that the rational curve on
$\tilde{X_0}$ which replaces the bad point lies on the inverse image of
this component, and meets the other components in a single point. Figure 2
sums this up much more clearly.

\begin{figure}
\begin{center}
{\tt    \setlength{\unitlength}{0.92pt}
\begin{picture}(348,123)
\thinlines    \put(194,13){$C_1\times C_2$}
              \put(290,13){$C_2^{(2)}$}
              \put(196,90){$C_1^{(2)}$}
              \put(290,71){\line(0,-1){36}} 
              \put(234,62){\line(1,0){36}}
              \put(258,71){\line(0,-1){36}}
              \put(234,90){\line(1,0){36}}
              \put(279,36){\framebox(36,36){}}
              \put(234,36){\framebox(36,36){}}
              \put(234,78){\framebox(36,36){}}
              \put(108,78){\vector(1,0){67}}
              \put(98,41){$C_1$}
              \put(22,42){$C_2$}
              \put(35,105){\line(1,-1){58}}
              \put(68,105){\line(-1,-1){58}}
\end{picture}}
\end{center}
\caption{The components of the symmetric square of a nodal curve and their
incidences.}
\end{figure}
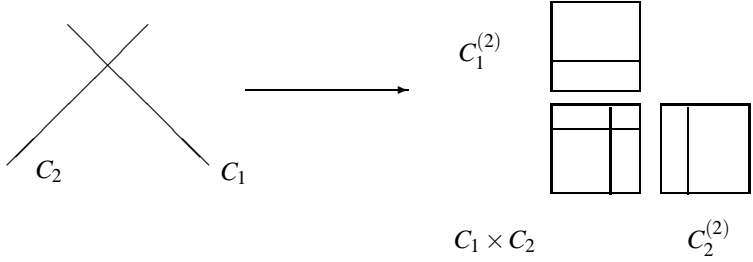

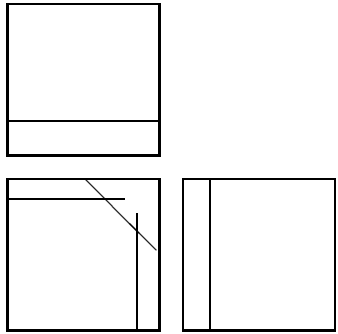
\begin{figure}
\begin{center}
{\tt    \setlength{\unitlength}{0.92pt}
\begin{picture}(154,155)
\thinlines    \put(93,72){\line(0,-1){62}}
              \put(10,96){\line(1,0){62}}
              \put(10,64){\line(1,0){48}}
              \put(63,58){\line(0,-1){48}}
              \put(42,72){\line(1,-1){29}}
              \put(82,10){\framebox(62,62){}}
              \put(10,10){\framebox(62,62){}}
              \put(10,82){\framebox(62,62){}}
\end{picture}}
\end{center}
\caption{Blowing up and regluing.}
\end{figure}

The product of $C$ with itself has components $C_i\times C_j$ for each pair
$C_i, C_j$ of irreducible components of $C$. The $\Z_2$-quotient turns each
product $C_i\times C_i$ into a component of $C^{(2)}$ isomorphic to the
symmetric square of $C_i$. The other components are not always, however,
products of curves. The component $C_i\times C_j$ ($i\neq j$) of $C^2$ is
identified with the component $C_j\times C_i$ by the $\Z_2$-action, but
because of incidences with other components, some points of the resulting
component are pinched together. The components of the normalization, however,
are not pinched, and are indeed isomorphic to products of curves blown up
at the points corresponding to bad points of $X_0$. 

The simplest example of this pinching is the symmetric square of a curve
$C=C_1\cup C_2$ with two components glued to each other at two points.
The symmetric square $C^{(2)}$ is Cohen-Macaulay, but the component
of $C^{(2)}$ which is not a symmetric square is not. However,
the normalization of $C^{(2)}$ is the disjoint union of $C_1^{(2)}$,
$C_2^{(2)}$ and $C_1\times C_2$. 

In the case of two irreducible components joined at a single node 
corresponding to a point $p$ on $C_1$ and $q$ on $C_2$, one may
describe the symmetric square of $C$ as $C_1\times C_2$ blown up in the point
$(p,q)$ with $C_1^{(2)}$ glued along the strict transform of $C_1\times\{q\}$
and $C_2^{(2)}$ glued along the strict transform of $\{p\}\times C_2$.

\section{Stability}

\subsection{Products and symmetric products of pointed curves}

In this section, the symbols $C$ or $C_i$ will denote smooth curves of
genus $g$ or $g_i$, and
$\delta$ or $\delta_i$ will denote reduced divisors on these curves. We denote
the degree of $\delta$ by $|\delta|$, since we have no occasion to refer to
the complete linear system usually denoted $|\delta|$.
It is
well known that the pair $(C,\delta)$ is of log-general type (that is, 
$K_C+\delta$ is big, which is the same as ample for curves) if 
$2g-2+|\delta|>0$. We define log surfaces related to these pointed curves.

\begin{definition}
The {\em log product} of $(C_1,\delta_1)$ with $(C_2,\delta_2)$ is the log
surface $(C_1\times C_2,\delta_1\times C_2 \cup C_1\times \delta_2)$. Since
we will never associate more than a single divisor to a single curve, the
notation $C_1\times_{\mrm{log}} C_2$ for the log product will not be 
confusing.

The {\em log symmetric square} of $(C,\delta)$, denoted $C^{(2)}_\mrm{log}$,
is defined as follows. Denote by $\fn{\pi}{C^2}{C^{(2)}}$ the quotient map for the
$\Z_2$-action swapping the factors. Then
$C^{(2)}_\mrm{log}$ is the surface 
$(C^2,\pi_*(\delta\times C\cup C\times\delta))$.
\end{definition}

Recall that the log-canonical class of a log variety $(X,D)$ is defined to
be $K_X+D$, where $K_X$ is the canonical class of $X$. $D$ is often called
a {\em boundary divisor}. It is clear that
the log-canonical class of $C_1\times_{\mrm{log}}C_2$ is ample if and only
if the log-canonical classes of $(C_1,\delta_1)$ and $(C_2,\delta_2)$ are
ample. The case of symmetric squares requires more work.

\begin{proposition}\label{clgt}
The log symmetric square of $(C,\delta)$ has ample log-canonical class exactly 
in the following situations:
\begin{enumerate}
\item $C$ is rational and $|\delta|>3$;
\item $C$ is elliptic and $|\delta|>2$;
\item $C$ is genus 2 and $|\delta|>1$;
\item $C$ is genus 3 hyperelliptic and $|\delta|>0$;
\item $C$ is non-hyperelliptic of genus 3;
\item $C$ is genus 4 or higher.
\end{enumerate}
The log-canonical class is nef if 
\begin{enumerate}
\item $C$ is rational and $|\delta|>2$;
\item $C$ is elliptic and $|\delta|>1$;
\item $C$ is genus 2 or higher.
\end{enumerate}
\end{proposition}

\begin{proof}
First, some notation is in order. Let $g$ be the genus of $C$. Denote by $B$ 
the boundary divisor of 
$C^{(2)}_\mrm{log}$. Denote its irreducible components by $B_i$, where
$i$ runs from 1 to $2|\delta|$. Let $\fn{\pi}{C^2}{C^{(2)}}$ be the quotient
map, and $\Delta$ the diagonal of $C^2$. Finally, let $D_i$ be the inverse
image of $B_i$ under $\pi$ for every $i$, and $D$ the union of the $D_i$. 
The $D_i$ are fibers of the projections from $C^2$ to its factors. 

Since $\pi$ is a finite map, and $\pi^*(K_{C^{(2)}}+B)=K_{C^2}+D-\Delta$, it
suffices to check that $K_{C^2}+D-\Delta$ is ample on $C^2$. For brevity,
we denote $K_{C^2}$ simply as $K$. We will use the Nakai-Moishezon criterion:
a divisor is ample if and only if its self-intersection is positive and it is 
positive on
every irreducible curve.
The following are easy to check:
\begin{eqnarray*}
K^2&=&2(2g-2)^2 \\
D^2&=&2|\delta|^2 \\
\Delta^2&=&2-2g \\
K.D&=&2|\delta|(2g-2) \\
K.\Delta&=& 4g-4 \\
D.\Delta&=&2|\delta|
\end{eqnarray*}
It follows that
\[
(K+D-\Delta)^2=2(2g-2)^2+(4|\delta|-5)(2g-2)+2|\delta|(|\delta|-2).
\]
This yields the following necessary conditions for ampleness:
\begin{enumerate}
\item If $g=0$, then $|\delta|>3$.
\item If $g=1$, then $|\delta|>2$.
\item If $g=2$, then $|\delta|>0$.
\end{enumerate}

First, we show that the conditions for $g=0$ or 1 are also sufficient.
If $g=0$, then $K+D-\Delta$ is a divisor of type $(|\delta|-3,|\delta|-3)$
on $\Pro^1\times\Pro^1$, hence ample as soon as $|\delta|>3$. If $g=1$, 
one has that a divisor on an abelian surface with positive self-intersection 
is ample or anti-ample (Corollary 2.2 of \cite{kani:ecas}). Since $K+D-\Delta$
is positive, for example, on $\Delta$, it is an ample divisor.

Suppose now that $C$ is hyperelliptic, and let $\Gamma$ be the graph of
the hyperelliptic involution in $C^2$. We have $\Gamma^2=2-2g$, and by the
Hurwitz formula
\[
(K+D-\Delta).\Gamma=-6+2g+2|\delta|
\]
which proves the necessity of all of the conditions given in the hypotheses.

If $C$ is a genus 2 curve, then $C^{(2)}$ is the blowup of an abelian
surface at a single point, and $\Gamma$ covers the exceptional curve for
this blowup. It follows that $K-\Delta$ is numerically equivalent to $\Gamma$.
Since $\Gamma$ is effective, $\Gamma+D$ is positive on any curve except
possibly its own components as soon as $|\delta|>0$, since then $D$ contains
a fiber in each direction. The computations above show that $\Gamma+D$ is
positive on $\Gamma$ once $|\delta|>1$, and for any component $D_i$ of $D$,
\[
(K+D-\Delta).D_i=2g-4+|\delta|
\]
by adjunction and geometric considerations. This finishes the genus 2 case.

Finally, it is well-known that the symmetric square of a genus 3 
nonhyperelliptic curve, or of a curve of genus 4 or higher already has an
ample canonical class, so no condition is necessary. Further, it is known
that the only (-2)-curve on $C^{(2)}$ for a genus 3 hyperelliptic curve is
the curve covered by $\Gamma$, so $\Gamma$ is the only curve to check for
positivity on $K+D-\Delta$, and this has been checked above.

The assertions about when the log-canonical class is nef follow from the
above computations.
\end{proof}

\subsection{Main theorems on stability}

Let $C$ be a nodal curve with components $C_i$. To each $C_i$, associate
a number $\delta_i$ which equals the number of nodes of $C_i$, with nodes
resulting from self-intersection counted twice. By genus of an irreducible
nodal curve, we mean the genus of the normalization, not the aritheoremetic
genus. The following lemma is necessary to reduce much of the proof of
the main theorem to the results of the previous section.

\begin{lemma}
The normalization of a product of irreducible nodal curves is the product of 
the normalizations of the factors. The normalization of the symmetric square
of an irreducible nodal curve is the symmetric square of the its normalization.
\end{lemma}

\begin{proof}
All
normalization morphisms will be denoted $\nu$. The normalization of a variety
$X$ will be denoted $X^\nu$. Let $C$, $C_1$, and $C_2$ be irreducible nodal
curves.

By the universal property of normalization, there exists a unique morphism
$\phi$ completing the diagram
\[
\xymatrix{
 & (C_1\times C_2)^\nu\ar[d]^\nu \\
C_1^\nu\times C_2^\nu\ar[r]\ar[ur]^\phi & C_1\times C_2
}
\]
commutatively. By Zariski's main theorem, the fibers of $\phi$ are connected. 
By commutativity of the diagram, any positive dimensional fiber of $\phi$
lies in a fiber of $C_1^\nu\times C_2^\nu\rightarrow C_1\times C_2$, which is 
impossible. We
conclude that $\phi$ is a homeomorphism, and hence an isomorphism, since
$(C_1\times C_2)^\nu$ is normal.

We now recycle the notation $\phi$ to denote the morphism $C^\nu\times C^\nu
\rightarrow (C\times C)^\nu$.
Via $\phi$, we have an action of $\Z_2$ on $(C\times C)^\nu$, and the 
normalization morphism is equivariant. Denote by $\pi$ the induced morphism
from $(C\times C)^\nu/\Z_2$ to $(C\times C)/\Z_2$. By an argument exactly like
the one just given, $(C\times C)^\nu/\Z_2$ is the normalization of
$(C\times C)/\Z_2$. It follows that $C^{(2)}$ is normalized by 
$(C^\nu)^{(2)}$.
\end{proof}

\begin{theorem}
Let $C$ be a nodal curve. Then $\mrm{Hilb}_2(C)$ is a stable surface if
and only if
\begin{enumerate}
\item for every genus 0 component $C_i$, $\delta_i>3$,
\item for every genus 1 component $C_i$, $\delta_i>2$,
\item for every genus 2 component $C_i$, $\delta_i>1$,
\item for every component $C_i$ with genus 3 hyperelliptic normalization,
$\delta_i>0$.
\end{enumerate}
\end{theorem}

\begin{proof}
We have already seen that $\mrm{Hilb}_2(C)$ has slc singularities, and
that it is Gorenstein. It remains to check that its dualizing sheaf
$\omega$ is ample.

$\omega$ is ample if and only if its restriction to every irreducible
component of the normalization of $\mrm{Hilb}_2(C)$ is ample. 
These components are products and symmetric squares of
curves or possibly blowups of these. The restriction of $\omega$ to a
component which is a product or symmetric product is exactly the log canonical
divisor on these surfaces considered in Proposition \ref{clgt}. Therefore it 
remains to
consider the blown-up components.

On the other hand, $K_X$ is positive on the exceptional curves of the 
blowups by the description of $\mrm{Hilb}_2(C)$ as the relative
canonical model of the Hilbert-Chow morphism.
\end{proof}

\begin{remark}
From Proposition \ref{clgt}, we can see that the morphism
$\mrm{Hilb}_S^2(C)\rightarrow S$ is relatively minimal unless there is
a genus 1 component $C_i$ with $\delta_i=1$. In this case, the symmetric
square $C_i^{(2)}$ is a ruled surface over $C_i$, and there is a unique
component $C_j$ of $C$ meeting $C_i$. Running the relative minimal model
program collapses $C_i^{(2)}$ to its curve of intersection with the blow
up $S$ of the product $C_i\times C_j$. In general, this will make the 
(-1)-curve on $S$ zero on $K_X$, so it will be blown down upon taking the
relative canonical model.
\end{remark}

From the preceeding theorem and remark, we have an explicit geometric
description of the relative minimal model for $\mrm{Hilb}_S^2(C)$ over
$S$. An explicit description of the canonical model is more complicated;
we are content to see a few examples.

\begin{example}
The simplest case to examine is when $C$ is a genus 3 hyperelliptic curve.
In this case, the graph of the hyperelliptic involution covers a rational
curve of self-intersection -2 on $C^{(2)}$. This is blown down to obtain
the canonical model. If $C$ occurs in a family of smooth hyperelliptic curves,
then these (-2)-curves sweep out a divisor which is collapsed by taking
the relative canonical model.
\end{example}

\begin{example}
If $C$ has a genus 2 component $C_i$ with $|\delta_i|=1$, then there
is a rational curve covered by the hyperelliptic involution (the unique
rational curve on $C_i^{(2)}$) with self-intersection -1. This
curve, however, is zero on $K_X$, so taking a minimal model does not
contract it. It is contracted to a smooth point by taking the canonical
model. 

For a concrete example, suppose $C$ is a genus 3 stable curve with two
components, a smooth genus 2 curve $C_1$ and an elliptic tail $C_2$. The
relative MMP collapses $C_2^{(2)}$ to its intersection with the blowup $S$
of $C_1\times C_2$. Now taking the relative canonical model blows down
the (-1)-curve on $S$ and the (-1)-curve on $C_1^{(2)}$. Thus
the stable limit of some smoothing of $C^{(2)}$ is the product $C_1\times C_2$
glued to an abelian surface (the Jacobian of $C_1$) along a curve isomorphic
to $C_2$. The canonical class is ample: its restriction to $C_1\times C_2$
is the tensor product of pullbacks of ample divisors from each factor, and
its restriction to $C_1^{(2)}$ is the class of a genus 2 curve, which has
positive self-intersection by the adjunction formula, hence is ample or
anti-ample. It is easy to see that the class of this curve is ample.
\end{example}

\begin{example}
In general, taking the canonical model introduces more singularities. If
$C$ has a rational component with only three special points, 
$\mrm{Hilb}_2(C)$ has a component isomorphic to $\Pro^2$, and such that
the canonical class restricted to this component is the canonical class of
$\Pro^2$ twisted by the coordinate axes (up to linear equivalence), hence
trivial. Therefore the morphism to the relative canonical model collapses
this component to a point, therefore also collapsing the curves of intersection
of this component with another component to points. The reader may find it
amusing to find the stable limit of a degeneration of smooth genus 3 curves
to a curve obtained by gluing a rational curve to an elliptic curve in three
distinct points.
\end{example}

\section{Global moduli results}

The study of this seemingly special degeneration gives all possible 
degenerations of symmetric squares. Indeed, let $X\rightarrow B$ be
a degeneration whose general fiber is the symmetric square of a curve
$C$. By the results of \cite{fan:dsp}, away from the special point
$0\in B$, $X$ comes from a family $C\rightarrow B\backslash\{0\}$ of
smooth curves. One then (possibly after base change) completes this family
to a family of stable curves $C\rightarrow B$. Then, the relative Hilbert
scheme may be formed, and the stable reduction process applied as above.

Let $g>2$.
Denote by $M_{g^{(2)}}$ the irreducible component of the moduli space of
stable surfaces containing the moduli point of the symmetric square of
a smooth genus $g$ curve. Let $M_g$ be the moduli space of genus $g$ stable
curves.

\begin{theorem}
There is a surjective, birational morphism $M_g\rightarrow M_{g^{(2)}}$ which
is an isomorphism over the locus of $M_g$ parameterizing smooth curves.
\end{theorem}

\begin{proof}
Let $H\rightarrow M_g$ be a surjective finite morphism from a smooth scheme 
$H$ induced by a family of stable curves. The existence of $H$ is guaranteed
by \cite{hm:mc}, Lemma 3.89. Let $C\rightarrow H$ be the family inducing
$H\rightarrow M_g$. Let $X\rightarrow H$ be relative Hilbert scheme of
length two subschemes of the fibers of $C\rightarrow H$. 

By \cite{karu:mmpbd}, Lemma 3.1, the relative canonical model 
$\tilde{X}\rightarrow H$ of $X\rightarrow H$ exists. By separatedness of
the moduli functor, the fibers of this relative canonical model coincide
with those obtained by taking relative canonical models of one-parameter
subfamilies.

$\tilde{X}$ induces a morphism $H\rightarrow M_{g^{(2)}}$ which  
descends to the desired morphism $\phi:M_g\rightarrow M_{g^{(2)}}$. The 
injectivity of $\phi$ on the locus of smooth curves in $M_g$ is the Torelli
theorem when $g=3$ and a theorem of Martens for higher genus curves 
\cite{cs:symprod}. The aforementioned result of \cite{fan:dsp} shows that
the image of the locus of smooth curves in $M_g$ is a dense open set in
$M_{g^{(2)}}$.
\end{proof}

Is $\phi$ also an isomorphism on the boundary? For a curve $C$, denote by
$C^{[2]}$ the special fiber of a relatively minimal model of the relative
Hilbert scheme for a smoothing as considered above. Denote the special 
fiber of the relative canonical model by $C^{\{2\}}$. Since minimal models
are isomorphic in codimension two, it follows that if 
$C^{\{2\}}\cong D^{\{2\}}$ for some curve $C$ and $D$, then 
$C^{[2]}\cong D^{[2]}$. Therefore we can avoid describing these relative
canonical models in our investigation of $\phi$.

\begin{lemma}
Suppose $C_i$ and $C_i'$ are smooth curves, $i=1, 2$. Then if
$C_1\times C_2\cong C_1'\times C_2'$, either
\begin{enumerate}
\item after possibly renumbering, $C_1\cong C_1'$ and $C_2\cong C_2'$, or
\item the $C_i$ and $C_i'$ are all elliptic curves. 
\end{enumerate}
\end{lemma}

\begin{proof}
First, it is clear that up to reordering, the genera of $C_i$ and $C_i'$ 
are equal.
By inclusion of a fiber, the isomorphism $C_1\times C_2\rightarrow 
C_1'\times C_2'$, and projection onto factors, we get finite covers
between the curves, for example $C_1\rightarrow C_1'$. Let $g$ be the common
genus of $C_1$ and $C_1'$. Then by Hurwitz,
\[
2g-2=n(2g-2)+r
\]
for positive integers $n$ and $r$. This is only possible if $2g-2=0$ or
if $n=1$ and $r=0$, from which the lemma follows.
\end{proof}

\begin{lemma}
Suppose $C$ and $D$ are smooth curves, and $C^{(2)}\cong D^{(2)}$. Then
either $C\cong D$ or $C$ and $D$ are genus 2.
\end{lemma}

\begin{proof}
See, e.g. \cite{cs:symprod}.
\end{proof}

\begin{lemma}
Suppose $C_1$, $C_2$, and $D$ are smooth curves.
Then $C_1\times C_2\not\cong D^{(2)}$. 
\end{lemma}

\begin{proof}
The invariants $\chi(\sh{O})$ and $K^2$ of a product of curves are related
by $K^2=8\chi$. This holds for symmetric squares only in the case $g=1$ where
$K^2=\chi=0$. The only product of curves with such invariants is a product
of elliptic curves, and the Kodaira dimension of a product of elliptic curves
differs from that of a symmetric square of an elliptic curve.
\end{proof}

\begin{remark}
Note, however, that it is possible that the blowup of a product of elliptic
curves be isomorphic to the symmetric square of a genus 2 curve.
\end{remark}

\begin{theorem}
The morphism $\phi$ defined above is an isomorphism.
\end{theorem}

\begin{proof}
Let $C$ be a stable curve. We reconstruct $C$ from $C^{[2]}$. For simplicity
assume that the components of $C$ are all smooth; the argument for 
self-intersecting components is similar but longer. Write $S$ for
$C^{[2]}$ in what follows. 

First of all, every pair of components of $S$ is either disjoint, meets in
a finite set of points, or meets along a curve. The only pairs of components 
that meet in a finite set are pairs of components isomorphic to symmetric
squares. Therefore, unless $C$ has an elliptic tail, we recover the 
symmetric squares of the components of $C$. Now if furthermore, $C$ has no
genus 2 components, these symmetric squares determine the components of
$C$. In any case, their incidences determine the incidences between 
components of $C$. 

In the case that $C$ has one or more genus 2 components, these may be
determined by the components of $S$ not isomorphic to a symmetric square,
since although $C_i^{(2)}$ may be isomorphic to $C_j^{(2)}$ for two distinct
genus 2 components, the component isomorphic to a blowup of $C_i\times C_j$
will determine $C_i$ and $C_j$, as well as their incidences with other
components of the curve.

A similar argument shows that the elliptic tails and their incidences can
be recovered, as long as the curve contains some component not of genus 1,
since a product $C_i\times C_j$ (the minimal model of some component of $S$)
will determine $C_i$ and $C_j$ as long as one of the two curves is not
elliptic.

It remains to cover the case of curves whose components are all elliptic 
curves. By the arguments above, the curve can be reconstructed, except
possibly for elliptic tails. Having found all of the components isomorphic
to symmetric squares, we look for components isomorphic to a product of
elliptic curves blown up at one point, which correspond to blown up products 
$T\times E$ of 
an elliptic tail $T$ with some other component $E$ of the curve. Call such
a component $S_1$. 
It is impossible this component to be the product of two elliptic tails 
(since such a curve would have genus 2). Therefore $S_1$ meets another
component $S_2$ which is a blow up of some product $T\times D$ of components
of $C$. $S_1$ is glued to $S_2$ along a curve isomorphic to $T$, so $T$
can be recovered. In this way all of the components of $C$ are recovered,
and their incidences can be deduced from the incidences of the components
of $S$.
\end{proof}

Therefore the moduli space of stable degenerations of a symmetric square of
a genus $g$ curve is simply $M_g$ as long as $g\geq 3$.

\end{document}